\input amstex

\documentstyle{amsppt}
  \magnification=1100
  \hsize=6.2truein
  \vsize=9.0truein
  \hoffset 0.1truein
  \parindent=2em

\NoBlackBoxes


\font\eusm=eusm10                   


\font\eusms=eusm7                       

\font\eusmss=eusm5                      


\newcount\theTime
\newcount\theHour
\newcount\theMinute
\newcount\theMinuteTens
\newcount\theScratch
\theTime=\number\time
\theHour=\theTime
\divide\theHour by 60
\theScratch=\theHour
\multiply\theScratch by 60
\theMinute=\theTime
\advance\theMinute by -\theScratch
\theMinuteTens=\theMinute
\divide\theMinuteTens by 10
\theScratch=\theMinuteTens
\multiply\theScratch by 10
\advance\theMinute by -\theScratch

\def\today{{\number\day\space
 \ifcase\month\or
  January\or February\or March\or April\or May\or June\or
  July\or August\or September\or October\or November\or December\fi
 \space\number\year}}

\define\Afr{{\frak A}}

\define\alphabar{{\overline\alpha}}

\define\Aut{\text{\rm Aut}}

\define\biggnm#1{
  \bigg|\bigg|#1\bigg|\bigg|}

\define\bignm#1{
  \big|\big|#1\big|\big|}

\define\Bt{{\widetilde B}}

\define\Cpx{\bold C}

\define\Dt{{\widetilde D}}

\define\eh{\hat e}

\define\Eto#1{E_{(\to{#1})}}

\define\fpamalg#1{{\dsize\;\operatornamewithlimits*_{#1}\;}}

\define\fpiamalg#1{{\tsize\;({*_{#1}})_{\raise-.5ex\hbox{$\ssize\iota\in I$}}}}

\define\freeprod#1#2{\mathchoice
     {\operatornamewithlimits{\ast}_{#1}^{#2}}
     {\raise.5ex\hbox{$\dsize\operatornamewithlimits{\ast}
      _{#1}^{#2}$}\,}
     {\text{oops!}}{\text{oops!}}}

\define\freeprodi{\mathchoice
     {\operatornamewithlimits{\ast}
      _{\iota\in I}}
     {\raise.5ex\hbox{$\dsize\operatornamewithlimits{\ast}
      _{\sssize\iota\in I}$}\,}
     {\text{oops!}}{\text{oops!}}}

\define\freeprodvni{\mathchoice
      {\operatornamewithlimits{\overline{\ast}}
       _{\iota\in I}}
      {\raise.5ex\hbox{$\dsize\operatornamewithlimits{\overline{\ast}}
       _{\sssize\iota\in I}$}\,}
      {\text{oops!}}{\text{oops!}}}

\define\Hil{{\mathchoice
     {\text{\eusm H}}
     {\text{\eusm H}}
     {\text{\eusms H}}
     {\text{\eusmss H}}}}


\define\Hilto#1{\Hil_{(\to{#1})}}

\define\id{\text{\rm id}}

\define\Integers{\bold Z}

\define\Ints{\Integers}

\define\Jup{^{(J)}}

\define\Keu{{\KHil}}

\define\KHil{{\mathchoice
     {\text{\eusm K}}
     {\text{\eusm K}}
     {\text{\eusms K}}
     {\text{\eusmss K}}}}

\define\ld#1{{\hbox{..}(#1)\hbox{..}}}

\define\Leu{{\mathchoice
     {\text{\eusm L}}
     {\text{\eusm L}}
     {\text{\eusms L}}
     {\text{\eusmss L}}}}

\define\lrnm#1{\left\|#1\right\|}

\define\lspan{\text{\rm span}@,@,@,}

\define\MvN{{\Cal M}}

\define\nm#1{\|#1\|}

\define\Nats{\Naturals}

\define\Naturals{{\bold N}}

\define\Oc{{\Cal O}}

\define\otdt{\otimes\cdots\otimes}

\define\otdts#1{\otimes_{#1}\cdots\otimes_{#1}}

\define\oup{^{\text{\rm o}}}

\define\owedge{{
     \operatorname{\raise.5ex\hbox{\text{$
     \ssize{\,\bigcirc\llap{$\ssize\wedge\,$}\,}$}}}}}

\define\owedgeo#1{{
     \underset{\raise.5ex\hbox
     {\text{$\ssize#1$}}}\to\owedge}}

\define\phibar{{\overline\phi}}

\define\Phit{{\widetilde\Phi}}

\define\psit{{\tilde\psi}}

\define\Psit{{\widetilde\Psi}}

\define\Pto#1{{P_{(\to{#1})}}}


\define\pup#1#2{{{\vphantom{#2}}^{#1}\!{#2}}\vphantom{#2}}

\define\QED{\newline
            \line{$\hfill$\qed}\enddemo}

\define\rank{\text{\rm rank}}

\define\restrict{\lower .3ex
     \hbox{\text{$|$}}}

\define\sigmat{{\widetilde\sigma}}

\define\smd#1#2{\underset{#2}\to{#1}}

\define\smdb#1#2{\undersetbrace{#2}\to{#1}}

\define\smdbp#1#2#3{\overset{#3}\to
     {\smd{#1}{#2}}}

\define\smdbpb#1#2#3{\oversetbrace{#3}\to
     {\smdb{#1}{#2}}}

\define\smdp#1#2#3{\overset{#3}\to
     {\smd{#1}{#2}}}

\define\smdpb#1#2#3{\oversetbrace{#3}\to
     {\smd{#1}{#2}}}

\define\smp#1#2{\overset{#2}\to
     {#1}}

\define\Tcirc{\bold T}

\define\Teu{\text{\eusm T}}

\define\tocdots
  {\leaders\hbox to 1em{\hss.\hss}\hfill}

\define\uh{{\hat u}}

\define\Veu{{\mathchoice
     {\text{\eusm V}}
     {\text{\eusm V}}
     {\text{\eusms V}}
     {\text{\eusmss V}}}}


  \newcount\mycitestyle \mycitestyle=1 

  \newcount\bibno \bibno=0
  \def\newbib#1{\advance\bibno by 1 \edef#1{\number\bibno}}
  \ifnum\mycitestyle=1 \def\cite#1{{\rm[\bf #1\rm]}} \fi
  \def\scite#1#2{{\rm[\bf #1\rm, #2]}}


  \newcount\ignorsec \ignorsec=0
  \def\notasec{\ignorsec=1}

  \newcount\secno \secno=0
  \def\newsec#1{\procno=0 \subsecno=0 \ignorsec=0
    \advance\secno by 1 \edef#1{\number\secno}
    \edef\currentsec{\number\secno}}

  \newcount\subsecno
  \def\newsubsec#1{\procno=0 \advance\subsecno by 1
    \edef\currentsec{\number\secno.\number\subsecno}
     \edef#1{\currentsec}}

  \newcount\appendixno \appendixno=0
  \def\newappendix#1{\procno=0 \ignorsec=0 \advance\appendixno by 1
    \ifnum\appendixno=1 \edef\appendixalpha{\hbox{A}}
      \else \ifnum\appendixno=2 \edef\appendixalpha{\hbox{B}} \fi
      \else \ifnum\appendixno=3 \edef\appendixalpha{\hbox{C}} \fi
      \else \ifnum\appendixno=4 \edef\appendixalpha{\hbox{D}} \fi
      \else \ifnum\appendixno=5 \edef\appendixalpha{\hbox{E}} \fi
      \else \ifnum\appendixno=6 \edef\appendixalpha{\hbox{F}} \fi
    \fi
    \edef#1{\appendixalpha}
    \edef\currentsec{\appendixalpha}}

  \newcount\procno \procno=0
  \def\newproc#1{\advance\procno by 1
   \ifnum\ignorsec=0 \edef#1{\currentsec.\number\procno}
                     \edef\currentproc{\currentsec.\number\procno}
   \else \edef#1{\number\procno}
         \edef\currentproc{\number\procno}
   \fi}

  \newcount\subprocno \subprocno=0
  \def\newsubproc#1{\advance\subprocno by 1
   \ifnum\subprocno=1 \edef#1{\currentproc a} \fi
   \ifnum\subprocno=2 \edef#1{\currentproc b} \fi
   \ifnum\subprocno=3 \edef#1{\currentproc c} \fi
   \ifnum\subprocno=4 \edef#1{\currentproc d} \fi
   \ifnum\subprocno=5 \edef#1{\currentproc e} \fi
   \ifnum\subprocno=6 \edef#1{\currentproc f} \fi
   \ifnum\subprocno=7 \edef#1{\currentproc g} \fi
   \ifnum\subprocno=8 \edef#1{\currentproc h} \fi
   \ifnum\subprocno=9 \edef#1{\currentproc i} \fi
   \ifnum\subprocno>9 \edef#1{TOO MANY SUBPROCS} \fi
  }

  \newcount\tagno \tagno=0
  \def\newtag#1{\advance\tagno by 1 \edef#1{\number\tagno}}



\notasec
 \newproc{\BigThm}
  \newtag{\htomegadelta}
  \newtag{\omegaalphaomega}
 \newproc{\freeprodfd}
 \newproc{\ZEFPdef}
 \newproc{\freepermfd}
 \newproc{\freeperminf}
 \newproc{\inclusionZEFP}
 \newproc{\exDpsi}
  \newtag{\statepsiOn}
 \newproc{\OinftyShift}
 \newproc{\CNTBV}
 \newproc{\CSZEFP}
 \newproc{\AllZEFP}
 \newproc{\AnyFinite}

\newbib{\BlanchardDykemaZZEmb}
\newbib{\BrownZZTopEntropy}
\newbib{\BrownChodaZZAppE}
\newbib{\ConnesNarnhoferThirring}
\newbib{\ConnesStormer}
\newbib{\CuntzZZOn}
\newbib{\DykemaZZExactFP}
\newbib{\KirchbergZZICM}
\newbib{\KirchbergPhillips}
\newbib{\SauvageotThouvenot}
\newbib{\StormerZZFreeShiftIIone}
\newbib{\StormerZZSurvey}
\newbib{\StormerZZShiftCstar}
\newbib{\VoiculescuZZSymmetries}
\newbib{\VoiculescuZZTopEntropy}
\newbib{\VDNbook}

\topmatter
  \title 
    Topological entropy of some automorphisms of reduced amalgamated free
    product C$^*$--algebras
  \endtitle

  \author Kenneth J\. Dykema \endauthor

  \date 20 May 1999 \enddate

  \rightheadtext{}
  \leftheadtext{}

  \address Dept.~of Mathematics,
           Texas A\&M University,
           College Station TX 77843-3368, USA
  \endaddress

  \email Ken.Dykema\@math.tamu.edu, {\it Internet URL:}
         http://www.math.tamu.edu/\~{\hskip0.1em}Ken.Dykema/
  \endemail

  \abstract
    Certain classes of automorphisms of reduced amalgamated free products of
    C$^*$--algebras are shown to have Brown--Voiculescu topological entropy
    zero.
    Also, for automorphisms of exact C$^*$--algebras, the
    Connes--Narnhofer--Thirring entropy is shown to be bounded above by the
    Brown--Voiculescu entropy.
    These facts are applied to generalize St\o{}rmer's result about the
    entropy of automorphisms of the II$_1$--factor of a free group.
  \endabstract


\endtopmatter


\document \TagsOnRight \baselineskip=18pt

\heading
\S1.  Introduction.
\endheading

Kolmogorov's entropy invariant was extended by Connes and
S\o{}rmer~\cite{\ConnesStormer} to an invariant $h_\tau(\alpha)$ for an
automorphism $\alpha$ of a von Neumann algebra with a given normal faithful
tracial state $\tau$ which is preserved by the automorphism.
One of the several results about the Connes--St\o{}rmer entropy
(see~\cite{\StormerZZSurvey} for a survey) is
St\o{}rmer's result~\cite{\StormerZZFreeShiftIIone} that the free shift on
$L(F_\infty)$ has entropy zero.
Here $L(F_\infty)$ is the II$_1$--factor defined by the left regular
representation of the free group $F_\infty$ on countably infinitely many
generators.
More generally, St\o{}rmer's theorem states that the entropy of $\sigma_*$ is
zero whenever $\sigma_*$ is the automorphism of $L(F_\infty)$ induced by a
permutation $\sigma$ of the generators of $F_\infty$ that has neither fixed
points nor finite cycles;
the free shift is the automorphism $\sigma_*$ where, when the generators of
$F_\infty$ are indexed by the integers, $\sigma$ corresponds to the shift
$n\mapsto n+1$.

The Connes--St\o{}rmer entropy was extended by Connes, Narnhofer and
Thirring~\cite{\ConnesNarnhoferThirring} to an invariant, generally referred to
as the CNT--entropy and denoted $h_\phi(\alpha)$, for an automorphism $\alpha$
of a unital C$^*$--algebra $A$ with respect to an $\alpha$--invariant state
$\phi$ of $A$.
See also Sauvageot's and Thouvenot's modification~\cite{\SauvageotThouvenot},
giving an entropy that is in general bounded above by the CNT--entropy and that
coincides with the CNT--entropy when the C$^*$--algebra $A$ is nuclear.
Theorem~VII.2 of~\cite{\ConnesNarnhoferThirring} shows that given
an automorphism $\alpha$ of a C$^*$--algebra $A$ preserving a state $\phi$, if
$\MvN$ is the von Neumann algebra generated by the image of $A$ under the GNS
representation of $\phi$, if $\alphabar$ and $\phibar$ are the canonical
extensions of $\alpha$ and  $\phi$ to $\MvN$, then
$h_\phibar(\alphabar)=h_\phi(\alpha)$.
(Their theorem is stated only for nuclear $A$ and hyperfinite $\MvN$, but their
proof applies generally.)

A noncommutative topological entropy was invented by
Voiculescu~\cite{\VoiculescuZZTopEntropy} for automorphisms of nuclear
C$^*$--algebras;
N\. Brown~\cite{\BrownZZTopEntropy} extended it to handle automorphisms of
exact C$^*$--algebras.
This Brown--Voiculescu entropy of an automorphism $\alpha$ is denoted
$ht(\alpha)$.
Voiculescu proved that if $\alpha$ is an automorphism of a unital nuclear
C$^*$--algebra $A$ and if $\phi$ is an $\alpha$--invariant state then
$h_\phi(\alpha)\le ht(\alpha)$.
Here we show (Proposition~\CNTBV) that the same inequality holds when $A$ is a
unital exact C$^*$--algebra.

In~\cite{\DykemaZZExactFP}, we proved that every reduced amalgamated free
product of exact C$^*$--algebras gives an exact C$^*$--algebra.
In this note, we build upon that proof to show that certain classes of
automorphisms of C$^*$--algebras arising as reduced amalgamated free products
have zero topological entropy.

The following section is the main part of the paper and contains the results
and their proofs.
At the end of it are two questions.

I would like to thank the members of the Institute of Mathematics in Luminy,
France and of the Erwin Schr\"odinger Insitute in Vienna, where much of this
research was done, for their hospitality;
I would like to thank especially J.B\. Cooper for organizing the Schr\"odinger
Institute's concentration in Fucntional Analysis in 1999.
Finally, the financial support of the CNRS of France and of the Schr\"odinger
Institute is gratefully acknowledged.

\heading
\S2.  Entropy of Automorphisms.
\endheading

\proclaim{Theorem \BigThm}
Let $B$ be a finite dimensional C$^*$--algebra, let $I$ be a set and for every
$\iota\in I$ let $A_\iota$ be a finite dimensional C$^*$--algebra containing
$B$ as a unital C$^*$--subalgebra and having a conditional expectation
$\phi_\iota:A_\iota\to B$ whose GNS representation is faithful.
Let
$$ (A,\phi)=\freeprodi(A_\iota,\phi_\iota) $$
be the reduced amalgamated free product of C$^*$--algebras and denote the
embeddings arising from the free product construction by
$\lambda_\iota:A_\iota\hookrightarrow A$.
Let $\sigma$ be a permutation of $I$ such that for every $\iota\in I$ there is
a $*$--isomorphism $\alpha_\iota:A_\iota\to A_{\sigma(\iota)}$ such that
$\alpha_\iota(B)=B$ and
$\phi_{\sigma(\iota)}\circ\alpha_\iota=\alpha_\iota\circ\phi_\iota$.
Assume further that the automorphism $\alpha_\iota{\restriction}_B$ of $B$ is
independent of $\iota\in I$, and call this automorphism $\beta$.
There is a unique automorphism $\alpha$ of $A$ such that
$\alpha\circ\lambda_\iota=\lambda_{\sigma(\iota)}\circ\alpha_\iota$ for all
$\iota\in I$.

Then $ht(\alpha)=0$.
\endproclaim
\demo{Proof}
In Voiculescu's construction~\cite{\VoiculescuZZSymmetries} of the reduced
amalgamated free product C$^*$--algebra $A$, one takes the Hilbert $B$--module
$E_\iota=L^2(A_\iota,\phi_\iota)$ on which $A_\iota$ acts via the GNS
representation, one lets $\xi_\iota=\widehat{1_{A_\iota}}\in E_\iota$, where
$A_\iota\ni a\mapsto\widehat a\in E_\iota$ is the defining map, one lets
$E_\iota\oup=E_\iota\ominus\xi_\iota B$, one constructs the free product of
Hilbert $B$--modules $(E,\xi)=\freeprodi(E_\iota,\xi_\iota)$, given by
$$ E=\xi B\oplus\bigoplus\Sb n\ge1\\
\iota_1,\iota_2\,\ldots,\iota_n\in I\\
\iota_1\ne\iota_2,\iota_2\ne\iota_3,\ldots,\iota_{n-1}\ne\iota_n \endSb
E_{\iota_1}\oup\otimes_BE_{\iota_2}\oup\otdts BE_{\iota_n}\oup, $$
and one defines $A$ acting on $E$;
(see~\scite{\DykemaZZExactFP}{\S1} for Voiculescu's construction in the
notation used here).
The $*$--isomorphism $\alpha_\iota:A_\iota\to A_{\sigma(\iota)}$ gives rise to
an invertible and isometric linear map $U_\iota:E_\iota\to E_{\sigma(\iota)}$
given by $U_\iota\widehat a=\widehat{\alpha_\iota(a)}$, (but note that
$U_\iota$ need not be $B$--linear).
Taking $A_\iota$, respectively $A_{\sigma(\iota)}$, acting via its GNS
representation on $E_\iota$, respectively $E_{\sigma(\iota)}$, we have for
$a\in A_\iota$ that $U_\iota aU_\iota^{-1}=\alpha_\iota(a)$.
Having assumed that $\alpha_\iota{\restriction}_B=\beta$ is independent of
$\iota$, we see that the collection of isometries $(U_\iota)_{\iota\in I}$
gives rise to an isometric and invertible linear map $U:E\to E$ given by
$U\xi b=\xi\beta(b)$ for $b\in B$ and
$U(\zeta_1\otdt\zeta_n)=(U_{\iota_1}\zeta_1)\otdt(U_{\iota_n}\zeta_n)$ for
$\zeta_j\in E_{\iota_j}\oup$ with $\iota_1,\ldots,\iota_n\in I$ and
$\iota_1\ne\iota_2,\ldots,\iota_{n-1}\ne\iota_n$.
The automorphism $\alpha$ of $A$ is then defined by $\alpha(x)=UxU^{-1}$.

Let $\pi$ denote the inclusion, arising from the free product construction, of
$A$ in $\Leu(E)$.
We will show that $ht(\pi,\alpha)=0$, and in order to do so we must show that
$ht(\pi,\alpha,\omega,\delta)=0$ for every finite subset $\omega$ of $A$ and
every $\delta>0$.
But for this it will suffice to let $\omega$ be a finite subset of any given
set whose linear span is a dense subset of $A$.
The set $W$ of reduced words in $(A_\iota)_{\iota\in I}$ has dense linear span
in $A$, and we will take $\omega\subseteq W$, where a {\it reduced word} is (an
element of $A$ given by) an expression of the form $a_1a_2\cdots a_n$, where
$n\ge1$, $a_j\in A_{\iota_j}\cap\ker\phi_{\iota_j}$ and
$\iota_1\ne\iota_2,\ldots,\iota_{n-1}\ne\iota_n$;
we call $n$ the {\it length} of the reduced word and we call the set
$\{\iota_1,\ldots,\iota_n\}\subseteq I$ the {\it alphabet} for the word;
we consider elements of $B$ to be reduced words of length $0$ and with alphabet
equal to the empty set.
If $\omega\subseteq W$ we define the alphabet for $\omega$ to be the union of
the alphabets of the elements of $\omega$.

Let $J$ be a subset of $I$ and let
$(A\Jup,\phi\Jup)=\freeprod{\iota\in J}{}(A_\iota,\phi_\iota)$ be the reduced
amalgamated free product of the subfamily.
Then $A\Jup$ acts canonically on the Hilbert $B$--module $E\Jup$, where
$(E\Jup,\xi)=\freeprod{\iota\in J}{}(E_\iota,\xi_\iota)$.
We will presently show in detail that $A\Jup$ is naturally embedded into $A$
and that there is a conditional expectation from $A$ onto $A\Jup$.
Note that $E\Jup$ is a complemented submodule of $E$;
let $\Theta\Jup:\Leu(E)\to\Leu(E\Jup)$ be given by compression.
Consider the Hilbert $B$--module
$$ E(J)=\eta B\oplus\bigoplus\bigoplus\Sb n\ge1\\
\iota_1,\iota_2\,\ldots,\iota_n\in I\\
\iota_1\ne\iota_2,\iota_2\ne\iota_3,\ldots,\iota_{n-1}\ne\iota_n\\
\iota_1\notin J \endSb
E_{\iota_1}\oup\otimes_BE_{\iota_2}\oup\otdts BE_{\iota_n}\oup, $$
where $\eta B$ is simply a copy of $B$ considered as a Hilbert $B$--module with
$\eta$ denoting the identity element of $B$.
There is then a unitary $V_J:E\to E\Jup\otimes_BE(J)$ given by erasing
parenthesis and absorbing $\eta$, analogous to the unitary
$E\to E_\iota\otimes_BE(\iota)$ in Voiculescu's construction of the reduced
amalgamated free product;
this unitary provides an embedding $i\Jup:\Leu(E\Jup)\to\Leu(E)$ given by
$i\Jup(x)=V_J^*(x\otimes1)V_J$, which then satisfies that
$\Theta\Jup\circ i\Jup$ is the identity on $\Leu(E\Jup)$.
Moreover, note that $i\Jup$ takes a reduced word considered as an element of
$A\Jup$ to the same reduced word considered as an element of $A$.
Hence $A\Jup$ is embedded into $A$ via $i\Jup$, and $\Theta\Jup$ provides a
conditional expectation from $A$ onto the embedded copy of $A\Jup$.

Let $\omega\subseteq W$ be a finite set of reduced words and let $\delta>0$;
we will find an upper bound for $rcp(\pi,\omega,\delta)$.
Let $q$ be the maximum of the lengths of the words belonging to $\omega$ and
let $J$ be the alphabet for $\omega$, which is thus a finite set.
Given $k\in\Nats$, consider the complemented submodule of $E\Jup$,
$$ \Eto k\Jup=\xi B\oplus\bigoplus\Sb 1\le n\le k\\
\iota_1,\iota_2\,\ldots,\iota_n\in J\\
\iota_1\ne\iota_2,\iota_2\ne\iota_3,\ldots,\iota_{n-1}\ne\iota_n \endSb
E_{\iota_1}\oup\otimes_BE_{\iota_2}\oup\otdts BE_{\iota_n}\oup $$
and let $\Phi\Jup_k:\Leu(E\Jup)\to\Leu(\Eto k\Jup)$ be given by compression.
In~\scite{\DykemaZZExactFP}{3.1}, unital completely positive maps
$\Psi_k\Jup:\Leu(\Eto k\Jup)\to\Leu(E\Jup)$ were constructed so that for every
$a\in A\Jup$, $\lim_{k\to\infty}\nm{a-\Psi_k\Jup\circ\Phi_k\Jup(a)}=0$.
Furthermore, from the proof of~\scite{\DykemaZZExactFP}{3.1} we see that for
every $\epsilon>0$ and every $q\in\Nats$ there is $k_0(\epsilon,q)\in\Nats$
such that for every reduced word $a\in A\Jup$ of length no greater than $q$, if
$k\ge k_0(\epsilon,q)$ then
$\nm{a-\Psi_k\Jup\circ\Phi_k\Jup(a)}\le\epsilon\nm a$;
moreover, $k_0(\epsilon,q)$ is universal, in the sense that it is the same for
all $J$ and all families $\bigl((A_\iota,\phi_\iota)\bigr)_{\iota\in J}$.
Hence, under the same conditions,
$\nm{a-i\Jup\circ\Psi_k\Jup\circ\Phi_k\Jup\circ\Theta\Jup(a)}\le\epsilon\nm a$.
Let us write $\Phit\Jup_k$ for the composition
$\Phi_k\Jup\circ\Theta\Jup:\Leu(E)\to\Leu(\Eto k\Jup)$ and
$\Psit\Jup_k$ for the composition
$i\Jup\circ\Psi_k\Jup:\Leu(\Eto k\Jup)\to\Leu(E)$.
Let $\epsilon=\delta/\max\{\nm a\mid a\in\omega\}$, let $q$ be the maximum of
the lengths of the words belonging to $\omega$ and let $k=k_0(\epsilon,q)$.
Since $J$ is a finite set and since each $E_\iota$ is a finite dimensional
complex vector space, the Hilbert $B$ module $\Eto k\Jup$ is a finite
dimensional vector space;
hence the C$^*$--algebra $\Leu(\Eto k\Jup)$ is finite dimensional.
Taking the unital completely positive maps $\Phit_k\Jup$ and $\Psit_k\Jup$, we
see that $rcp(\pi,\omega,\delta)\le\rank\bigl(\Leu(\Eto k\Jup)\bigr)$.
We now perform a crude (but sufficient) estimate of this rank.
Let $d(J)$ be the maximum over $\iota\in J$ of the dimension of $E_\iota$ as a
vector space;
then we can estimate
$$ \dim(\Eto k\Jup)\le\dim(B)+\sum_{n=1}^k|J|^nd(J)^n
\le\dim(B)+k|J|^kd(J)^k. $$
Let $\rho$ be a faithful representation of $B$ on a finite dimensional Hilbert
space $\Veu$.
Then the C$^*$--algebra $\Leu(\Eto k\Jup)$ is faithfully represented on the
Hilbert space $\Eto k\Jup\otimes_\rho\Veu$, which has dimension
$\le\dim(\Eto k\Jup)\dim(\Veu)$.
Thus we have
$$ rcp(\pi,\omega,\delta)\le\bigl(\dim(B)+k|J|^kd(J)^k\bigr)\dim(\Veu). $$

Now we are in a position to show that
$$ ht(\pi,\alpha,\omega,\delta)=0. \tag{\htomegadelta} $$
Given the nature of our automorphism $\alpha$, for every $n\in\Nats$ the
maximum length and the maximum norm of words belonging to
$$ \omega\cup\alpha(\omega)\cup\cdots\cup\alpha^{n-1}(\omega)
\tag{\omegaalphaomega} $$
are the same as for $\omega$, and we may choose $k=k_0(q,\epsilon)$ as for
$\omega$ above.
However, the alphabet $J_n$ of the set of words~(\omegaalphaomega) is equal to
$J\cup\sigma(J)\cup\cdots\cup\sigma^{n-1}(J)$, and thus $|J_n|\le n|J|$.
But the existence of the isomorphisms $\alpha_\iota$ preserving
conditional expectations implies that $\dim(E_{\sigma(\iota)})=\dim(E_\iota)$,
and hence $d(J_n)=d(J)$.
Hence we have the estimate
$$ rcp(\pi,\omega\cup\alpha(\omega)\cup\cdots\cup\alpha^{n-1}(\omega),\delta)
\le\bigl(\dim(B)+kn^k|J|^kd(J)^k\bigr)\dim(\Veu). $$
As the upper bound grows subexponentially in $n$, the estimate
implies~(\htomegadelta).
\QED

We now list as corollaries some particular sorts of automorphisms to which the
above theorem applies.
First we have free products of automorphisms, which correspond to when the
permutation $\sigma$ in Theorem~\BigThm{} is the identity.

\proclaim{Corollary \freeprodfd}
Let
$$ (A,\phi)=\freeprodi(A_\iota,\phi_\iota) $$
be the reduced amalgamated free product of finite dimensional C$^*$--algebras
as in the statement of Theorem~\BigThm.
For every $\iota\in I$ let $\alpha_\iota\in\Aut(A_\iota)$ be such that
$\alpha_\iota(B)=B$,
$\phi_\iota\circ\alpha_\iota=\alpha_\iota\circ\phi_\iota$;
suppose that the automorphism $\alpha_\iota{\restriction}_B$ of $B$ is the same
for all $\iota\in I$.
Let $\alpha=\freeprodi\alpha_\iota\in\Aut(A)$;
by this we mean that $\alpha$ is the automorphism of $A$ that when restricted
to the naturally embedded copy of $A_\iota$ in $A$ is $\alpha_\iota$.

Then $ht(\alpha)=0$.
\endproclaim

Next we have the free shifts and their analogues for general permutations.
\proclaim{Definition \ZEFPdef}\rm
If $(A,\phi)=\freeprodi(A_\iota,\phi_\iota)$ is a reduced amalgamated free
product of C$^*$--algebras, where each $(A_\iota,\phi_\iota)$ is a copy of a
fixed pair $(D,\psi)$ of a unital exact C$^*$--algebra $D$ and a conditional
expectation $\psi$ from $D$ onto a unital C$^*$--subalgebra $B$ having faithful
GNS representation, and if $\sigma$ is a permutation of the index set $I$, then
what we call the corresponding {\it free permutation} is the automorphism
$\sigma_*$ of $A$ sending the embedded copy of $A_\iota$ in $A$ identically
to the embedded copy of $A_{\sigma(\iota)}$ in $A$, for every $\iota\in I$.

We say that the pair $(D,\psi)$ has the {\it ZEFP property} (with respect to
$ht$) if $ht(\sigma_*)=0$ whenever $\sigma_*$ is a free permutation of a free
product of some copies of $(D,\psi)$.
\endproclaim
The acronym ZEFP is for ``zero entropy free permutation.''

\proclaim{Corollary \freepermfd}
Let $B$ and $D$ be finite dimensional C$^*$--algebras with $B$ contained as a
unital C$^*$--subalgebra of $D$;
let $\psi:D\to B$ be a conditional expectation whose GNS representation is
faithful.
Then $(D,\psi)$ has the ZEFP property.
\endproclaim

\proclaim{Corollary \freeperminf}
Let $J$ be a set, let $B$ be a finite dimensional C$^*$--algebra and for every
$\iota\in J$ let $D_\iota$ be a finite
dimensional C$^*$--algebra and $\psi_\iota:D_\iota\to B$ is a conditional
expectation having faithful GNS representation.
Let  $(D,\psi)=\freeprod{\iota\in J}{}(D_\iota,\psi_\iota)$.
Then $(D,\psi)$ has the ZEFP property
\endproclaim
\demo{Proof}
If $I$ is a set and if $\sigma$ is a permutation of $I$, let $\sigma_*$ be the
corresponding free permutation of the free product of $|I|$ copies of
$(D,\psi)$.
Then $\sigma_*$ is in the obvious way equal to a free permutation of a reduced
free product of finite dimensional C$^*$--algebras, corresponding to the
permutation $\sigma\times\id_J$ of $I\times J$.
Thus $ht(\sigma_*)=0$ by Theorem~\BigThm.
\enddemo

\proclaim{Definition and Proposition \inclusionZEFP}
Let $(D,\psi)$ and $(\Dt,\psit)$ be pairs of a unital exact C$^*$--algebras $D$
and $\Dt$ with conditional expectations $\psi$ from $D$ onto a unital
C$^*$--subalgebra $B\subseteq D$ and $\psit$ from $\Dt$ onto a unital
C$^*$--subalgebra $\Bt\subseteq\Dt$, whose GNS representations are faithful.
We say $(D,\psi)$ is {\rm included} in $(\Dt,\psit)$, and write
$(D,\psi)\subseteq(\Dt,\psit)$, if $D$ is a C$^*$--subalgebra of $\Dt$ in such
a way that $B\subseteq\Bt$ and $\psit{\restriction}_D=\psi$.
We call the inclusion $(D,\psi)\subseteq(\Dt,\psit)$ {\rm unital} if $D$ is a
unital C$^*$--subalgebra of $\Dt$.

If $(D,\psi)\subseteq(\Dt,\psit)$ and if $(\Dt,\psit)$ has the ZEFP property
then $(D,\psi)$ has the ZEFP property.
\endproclaim
\demo{Proof}
First suppose that the inclusion is unital.
By the main result of~\cite{\BlanchardDykemaZZEmb}, the free product of $|I|$
copies of $(D,\psi)$ embeds in the free product of $|I|$ copies of
$(\Dt,\psit)$.
Let $\sigma$ be a permutation of $I$,
let $\sigma_*$ be corresponding free permutation of the free product of $|I|$
copies of $(D,\psi)$ and let $\sigmat_*$ be the free permutation of the free
product of $|I|$ copies of $(\Dt,\psit)$.
Then $\sigma_*$ is the restriction of $\sigmat_*$.
As the Brown--Voiculescu topological entropy is
monotone~\scite{\BrownZZTopEntropy}{2.1}, we have $ht(\sigma_*)=0$;
hence $(D,\psi)$ has the ZEFP property.

If the inclusion $(D,\psi)\subseteq(\Dt,\psit)$ is nonunital, let $p\in\Dt$
denote the identity element of $D$ and let $1$ denote the identity element of
$\Dt$;
then $1-p\in\Bt$.
Let $D'=D+\Cpx(1-p)\subseteq\Dt$ and let $B'=B+\Cpx(1-p)\subseteq\Bt$;
then for $d\in D$ and $\lambda\in\Cpx$,
$\psit(d+\lambda(1-p))=\psi(d)+\lambda(1-p)$;
let $\psi'=\psit{\restriction}_{D'}:D'\to B'$.
Then by the unital case just proved, $(D',\psi')$ has the ZEFP property.
Let $I$ be a set and let $(A',\phi')=\freeprodi(A'_\iota,\phi'_\iota)$ where
each $(A'_\iota,\phi'_\iota)$ is a copy of $(D',\psi')$;
let $(A,\phi)=\freeprodi(A_\iota,\phi_\iota)$ where each $(A_\iota,\phi_\iota)$
is a copy of $(D,\psi)$.
Then $p\in B'\in A'$ and $A$ is canonically isomorphic to $pA'p$;
if $\sigma_*$ is a free 
permutation on $A$ corresponding to a permutation $\sigma$ of $I$,
then $\sigma_*$ is the restriction of the corresponding free permutation
$\sigma'_*$ of $A'$ to $pA'p$.
Again by monotonicity, we see that $ht(\sigma_*)=0$ and $(D,\psi)$ has the ZEFP
property.
\QED

Application of Corollary~\freeperminf{} and Proposition~\inclusionZEFP{} leads
to many examples, a few of which are below.

\proclaim{Examples \exDpsi}
The following pairs have the ZEFP property.
\roster
\item"(i)" $(\Teu,\phi_1)$ where $\Teu$ is the Toeplitz algebra, which is
generated by a nonunitary isometry $v$, and where $\phi_1$ is the state on
$\Teu$ satisfying $\phi_1(vv^*)=0$;
\item"(ii)" $(\Oc_\infty,\phi)$ where $\Oc_\infty$ is the Cuntz
algebra~\cite{\CuntzZZOn}, which is generated by isometries $s_1,s_2,\ldots$
having orthogonal ranges, and where $\phi$ is the state on $\Oc_\infty$
such that $\phi(s_js_j^*)=0$ for all $j$;
\item"(iii)" $(\Oc_n,\phi_n)$, with $n\in\Nats$, $n\ge2$, where $\Oc_n$ is the
Cuntz algebra~\cite{\CuntzZZOn}, which is generated by isometries
$s_1,\ldots,s_n$, whose range projections sum to $1$, and where, for any fixed
choice of $\gamma_1,\ldots,\gamma_n\in[0,1]$ such that
$\gamma_1+\cdots\gamma_n=1$, $\phi_n$ is the state on $\Oc_n$ given by
$$ \phi_n(s_{i_1}s_{i_2}\cdots s_{i_k}s_{j_\ell}^*\cdots s_{j_2}^*s_{j_1}^*)=
\cases \gamma_{i_1}\gamma_{i_2}\cdots\gamma_{i_k}
&\text{if }k=\ell,\,i_1=j_1,\ldots,i_k=j_k \\
0&\text{otherwise;}\endcases \tag{\statepsiOn} $$
\item"(iv)" $(\Oc_\infty,\phi_\infty)$ where $\Oc_\infty$ is generated by
isometries $s_1,s_2,\ldots$ having orthogonal ranges and where for any fixed
choice of $\gamma_1,\gamma_2,\ldots\in[0,1]$ such that
$\sum_1^\infty\gamma_j\le1$, $\phi_\infty$ is the state on $\Oc_\infty$
satisfying~(\statepsiOn);
\item"(v)" $(C(\Tcirc),\tau)$ where $\Tcirc$ is the circle and where the state
$\tau$ is given by Lebesgue measure on $\Tcirc$.
\endroster
\endproclaim
\demo{Proof}
For~(i), let $D_1=\Cpx\oplus\Cpx$ with minimal projection $p\in D_1$ and let
$\psi_1$ be the state on $D_1$ such that $\psi_1(p)=1/2$;
let $D_2=M_2(\Cpx)$ with a system of matrix units $(e_{ij})_{1\le i,j\le 2}$ in
$D_2$ and let $\psi_2$ be the state on $D_2$ so that $\psi_2(e_{11})=1$.
Let $(\Dt,\psit)=(D_1,\psi_1)*(D_2,\psi_2)$.
Considering the unitary $u=1-2p\in D_1$, we see that $L^2(D_1,\psi_1)$ has
orthonormal basis $\{\widehat{1_{D_1}},\uh\}$;
moreover, $L^2(D_2,\psi_2)$ has orthonormal basis
$\{\widehat{1_{D_2}},\eh_{21}\}$.
Therefore, $L^2(\Dt,\psit)$ has orthonormal basis
$$ \align
\{\xi\}&\cup\{\uh,\,\uh\otimes\eh_{21},\,\uh\otimes\eh_{21}\otimes\uh,\,
\uh\otimes\eh_{21}\otimes\uh\otimes\eh_{21},\,\ldots\}\quad\cup \\
&\cup\quad
\{\eh_{21},\,\eh_{21}\otimes\uh,\,\eh_{21}\otimes\uh\otimes\eh_{21},\,
\eh_{21}\otimes\uh\otimes\eh_{21}\otimes\uh,\,\ldots\},
\endalign $$
where $\xi=\widehat{1_\Dt}$;
moreover, $\psit$ is the vector state associated to $\xi$.
Let $v=e_{21}ue_{22}+e_{11}ue_{21}\in\Dt$.
Then $v$ is an isometry satisfying
$$ \align
v:\,&\xi\mapsto\uh\otimes\eh_{21} \\
&\uh\otimes(\cdots)
 \mapsto\uh\otimes\eh_{21}\otimes\uh\otimes(\cdots) \\
&\eh_{21}\otimes(\cdots)
 \mapsto\eh_{21}\otimes\uh\otimes\eh_{21}\otimes(\cdots).
\endalign $$
Thus the C$^*$--subalgebra of $\Dt$ generated by $v$ is isomorphic to $\Teu$
and, as $\xi$ is orthogonal to the range space of $v$, the restriction of
$\psit$ to the copy of $\Teu$ is the state $\phi_1$ described in~(i).
Now Corollary~\freeperminf{} and Proposition~\inclusionZEFP{} imply that
$(\Teu,\phi_1)$ has the ZEFP property.

Note that~(ii) is a special case of~(iv).
However, for future reference we would like to point out how~(ii) follows
from~(i).
{}From~\scite{\VoiculescuZZSymmetries}{\S2} (or see~\scite{\VDNbook}{1.5.10}),
$(\Oc_\infty,\phi)$ is the free product of countably infinitely many
copies of $(\Teu,\phi_1)$.
Hence by Corollary~\freeperminf{} and Proposition~\inclusionZEFP,
$(\Oc_\infty,\phi)$ has the ZEFP property.

For~(iii), let $\Bt=\Cpx\oplus\Cpx$ with minimal projection $p$;
let $D_1=M_2(\Cpx)$ with a system of matrix units $(e_{ij})_{0\le i,j\le1}$,
with $\Bt$ unitally embedded by identifying $p$ and $e_{11}$, and with
conditional expectation $\psi_1:D_1\to\Bt$ given by
$$ \psi_1\bigl(\sum_{i,j=0}^1c_{ij}e_{ij}\bigr)=c_{11}p+c_{00}(1-p); $$
let $D_2=M_{n+1}(\Cpx)$ with a system of matrix units
$(f_{ij})_{0\le i,j\le n}$, with $\Bt$ unitally embedded by identifying $1-p$
and $f_{00}$ and with conditional expectation $\psi_2:D_2\to\Bt$ given by
$$ \psi_2\bigl(\sum_{i,j=0}^nc_{ij}f_{ij}\bigr)=
\bigl(\sum_{j=1}^n\gamma_jc_{jj}\bigr)p+c_{00}(1-p). $$
Let
$$ (\Dt,\psit)=(D_1,\psi_1)*(D_2,\psi_2). $$
For every $k\in\{1,\ldots,n\}$, let $s_k=f_{k0}e_{01}\in\Dt$.
Then $s_k^*s_k=p$ and $s_ks_k^*=f_{kk}$.
In $p\Dt p$, $s_1,\ldots,s_n$ are isometries with range projections
summing to $p$, so they generate a copy of $\Oc_n$ in $\Dt$ with identity
element $p$ and to which the conditional expectation $\psit$ restricts to a
state, $\phi_n$ (when $\Cpx p$ is identified with $\Cpx$).
It is clear that $\phi_n(s_js_j^*)=\gamma_j$;
in order to see that~(\statepsiOn) holds, one can argue by induction on $k$ and
use freeness.
Now Corollary~\freeperminf{} and Proposition~\inclusionZEFP{} imply that
$(\Oc_n,\phi_n)$ has the ZEFP property.

The proof of~(iv) is similar to the that of~(iii), but taking $D_2$ to be
the unitization of the C$^*$--algebra, $\Keu$, of compact operators on
separable infinite dimensional Hilbert space.
Letting $(f_{ij})_{i,j\ge0}$ be a system of matrix units for $\Keu$, embed
$\Bt$ in $D_2$ by identifying $1-p$ and $f_{00}$, and let $\psi_2:D_2\to\Bt$ be
the conditional expectation given by
$$ \psi_2(1)=1\qquad\psi_2(f_{jj})=\gamma_j p\quad(j\ge1)\qquad
\psi_2(f_{00})=1-p. $$
Then letting $s_j=f_{j0}e_{01}$, ($j\ge1$) we have $s_j^*s_j=p$ and
$s_js_j^*=f_{jj}$; hence $\{s_1,s_2\ldots\}$ generates a copy of $\Oc_\infty$
in $p\Dt p$, to which the restriction of $\psit$ is seen to be $\phi_\infty$ as
described in~(iv) above.

For~(v), it is only required to apply Corollary~\freeperminf{} and
Proposition~\inclusionZEFP{} after noting that the choice of an infinite order
element in the group $\Ints_2*\Ints_2$, (the free product of the
two--element group with itself), gives rise to an canonical trace preserving
embedding of the reduced group C$^*$--algebra
$C^*_r(\Ints)\cong C(\Tcirc)$ in the reduced group C$^*$--algebra
$C^*_r(\Ints_2*\Ints_2)$, which in turn arises as the reduced free product
$$ (C^*_r(\Ints_2*\Ints_2),\tau_{\Ints_2*\Ints_2})
=(C^*_r(\Ints_2),\tau_{\Ints_2})*(C^*_r(\Ints_2),\tau_{\Ints_2}) $$
of finite dimensional C$^*$--algebras.
\QED

Example~\exDpsi(i), can be used to give another proof of Brown's and Choda's
result that the free permutations on the Cuntz algebra $\Oc_\infty$ have
topological entropy zero.
\proclaim{Proposition~\OinftyShift}{\rm\;(\cite{\BrownChodaZZAppE})}
Let $\{s_1,s_2,\ldots\}$ be a family of isometries having
orthogonal ranges and generating the Cuntz algebra $\Oc_\infty$;
let $\alpha$ be an automorphism of $\Oc_\infty$ given by
$\alpha(s_k)=s_{\sigma(k)}$, where $\sigma$ is some permutation of $\Nats$.
Then $ht(\alpha)=0$.
\endproclaim
\demo{Proof}
As mentioned in the proof of~\exDpsi(ii) above, $\Oc_\infty$ is the free
product of countably infinitely many copies of $(\Teu,\phi_1)$ as
in~\exDpsi(i), indexed by $\Nats$.
We get $ht(\alpha)=0$ because $(\Teu,\phi_1)$ has the ZEFP property.
\QED

We will use Example~\exDpsi(v) to generalize, to the case of arbitrary
permutations, St\o{}rmer's result~\cite{\StormerZZFreeShiftIIone} about free
shifts on $L(F_\infty)$.
For this, we need to extend Voiculescu's inequality
$h_\sigma(\alpha)\le ht(\alpha)$ to the case of automorphisms of unital exact
C$^*$--algebras.
The proof below is inspired by
Voiculescu's~\scite{\VoiculescuZZTopEntropy}{4.6}; we refer
to~\cite{\ConnesNarnhoferThirring} and~\cite{\BrownZZTopEntropy} for relevant
concepts and definitions.

\proclaim{Proposition \CNTBV}
Let $A$ be a unital exact C$^*$--algebra, let $\alpha\in\Aut(A)$ and let
$\sigma$ be a state on $A$ satisfying $\sigma\circ\alpha=\sigma$.
Then $h_\sigma(\alpha)\le ht(\alpha)$.
\endproclaim
\demo{Proof}
Let $\gamma:M_k(\Cpx)\to A$ be a unital completely positive map.
Let $\omega$ be a finite subset of $A$ such that
$\gamma(M_k(\Cpx))\subseteq\lspan\omega$ and
$$ \gamma\bigl(\{x\in M_k(\Cpx)\mid\nm x\le1\}\bigr)\subseteq
\Bigl\{\sum_{x\in\omega}\lambda(x)x\;\Big|\;\lambda(x)\in\Cpx,\,
\sum_{x\in\omega}|\lambda(x)|\le1\Bigr\}; $$
for future reference, assume that also the identity element of $A$ belongs to
$\omega$.
Let $\pi:A\to\Leu(\Hil)$ be a faithful representation of $A$ on a Hilbert space
$\Hil$.
Let $\delta>0$ and $n\in\Nats$ and suppose that $D$ is a finite dimensional
C$^*$--algebra and that $\phi:A\to D$ and $\psi:D\to\Leu(\Hil)$ are unital
completely positive maps such that
$\forall a\in\omega\cup\alpha(\omega)\cup\cdots\cup\alpha^{n-1}(\omega)$,
$\nm{\psi\circ\phi(a)-\pi(a)}<\delta$.
Then for all $x\in M_k(\Cpx)$ with $\nm x\le1$ and for all
$j\in\{0,1,\ldots,n-1\}$,
$$ \nm{\psi\circ\phi\circ\alpha^j\circ\gamma(x)
-\pi\circ\alpha^j\circ\gamma(x)}<\delta. $$
Let $C$ be the C$^*$--algebra generated by $\pi(A)\cup\psi(D)$.
Consider an abelian model, call it $\Afr$, for
$\bigl(A,\phi,(\alpha^j\circ\gamma)_{j=0}^{n-1}\bigr)$ consisting
of an abelian finite dimensional C$^*$--algebra $B$, a unital completely
positive map $P:A\to B$, a state $\mu$ on $B$ such that $\mu\circ P=\sigma$
and $*$--subalgebras $B_1,\ldots,B_n$ of $B$.
There is a unital completely positive map $P':C\to B$ such that
$P'\circ\pi=P$.
If $E_j:B\to B_j$ are the canonical conditional expectations with respect to
$\mu$, then letting
$$ \align
\rho_j&=E_j\circ P\circ\alpha^j\circ\gamma:A\to B_j \\
\rho_j'&=E_j\circ P'\circ\psi\circ\phi\circ\alpha^j\circ\gamma:A\to B_j,
\endalign $$
we have $\nm{\rho_j-\rho_j'}\le\delta$ for all $j$.
Then by~\scite{\ConnesNarnhoferThirring}{IV.2},
$|s_\mu(\rho_j)-s_\mu(\rho_j')|<\eta$ where
$\eta=3\delta+6\delta\log(1+k^2\delta^{-1})$.
Let $\sigma'=\mu\circ P'$ and let $\Afr'$ be the abelian model for
$\bigl(C,\sigma',(\psi\circ\phi\circ\alpha^j\circ\gamma)_{j=0}^{n-1}\bigr)$
consisting of $(B,\mu,B_1,\ldots,B_n)$ and the completely positive map
$P':C\to B$.
Then from equation~(III.3) of~\cite{\ConnesNarnhoferThirring}, the entropy of
the abelian model $\Afr$ differs from that of $\Afr'$ by no more than $n\eta$.
Moreover, the entropy of the abelian model $\Afr'$ is bounded above by
$H_{\sigma'}\bigl((\psi\circ\phi\circ\alpha^j\circ\gamma)_{j=0}^{n-1}\bigr)$;
this is by~\scite{\ConnesNarnhoferThirring}{III.6(a,c)} bounded above by
$H_{\sigma'}(\psi)$, which is $\le\log\rank(D)$.
We may choose $(D,\phi,\psi)$ so that
$\rank(D)\le rcp(\pi,\omega\cup\alpha(\omega)\cup\cdots\cup\alpha^{n-1}(\omega),4\delta)$;
indeed, had we not required $\phi$ and $\psi$ to be unital, we could have
chosen $(D,\phi,\psi)$ so that
$\rank(D)=rcp(\pi,\omega\cup\alpha(\omega)\cup\cdots\cup 
\alpha^{n-1}(\omega),\delta)$,
but as $1\in\omega$,
any nonunital $\phi$ and $\psi$ can be rescaled to give unital ones.
Hence we find
$$ H_\sigma(\gamma,\alpha\circ\gamma,\cdots,\alpha^{n-1}\circ\gamma)\le
\log rcp(\pi,\omega\cup\cdots\cup\alpha^{n-1}(\omega),4\delta)+n\eta; $$
therefore $h_{\sigma,\alpha}(\gamma)\le ht(\pi,\alpha,\omega,4\delta)+\eta$.
If $\delta\to0$ then $\eta\to0$ and we find
$h_{\sigma,\alpha}(\gamma)\le ht(\alpha)$;
hence $h_\sigma(\alpha)\le ht(\alpha)$.
\QED

\proclaim{Corollary~\CSZEFP}
Let $\sigma_*$ be the automorphism of the II$_1$--factor $L(F_\infty)$ induced
by an arbitrary permutation $\sigma$ of the generators of the group $F_\infty$.
Then the Connes--St\o{}rmer entropy of $\sigma_*$ is zero.
\endproclaim
\demo{Proof}
Let $\tau$ be the tracial state on $L(F_\infty)$.
Combining Example~\exDpsi(v) with Proposition~\CNTBV, we find that the
CNT--entropy $h_\tau(\sigma_{r,*})$ is zero, where $\sigma_{r,*}$ is the
automorphism of $C^*_r(F_\infty)$ arising from the permutation $\sigma$ of the
generators of $F_\infty$ and where $\tau$ is the unique tracial state on
$C^*_r(F_\infty)$.
But $h_\tau(\sigma_{r,*})$ is equal to the CNT--entropy (hence, to the
Connes--St\o{}rmer entropy) of the corresponding automorphism $\sigma_*$ of
$L(F_\infty)$.
\QED

The following question is quite natural.
\proclaim{Question \AllZEFP}\rm
Does every pair $(D,\psi)$, where $D$ is a unital exact C$^*$--algebra and
where $\psi$ is a conditional expectation from $D$ onto a unital
C$^*$--subalgebra, have the ZEFP property?
\endproclaim

This seems like an appropriate place to point out that by recent work of
Kirchberg~\cite{\KirchbergZZICM}, \cite{\KirchbergPhillips},  with $(D,\psi)$
as in Question~\AllZEFP, one can always realize $D\subseteq\Oc_2$;
if one could realize $(D,\psi)\subseteq(\Oc_2,\phi_2)$, with $(\Oc_2,\phi_2)$
as in Example~\exDpsi(iii), then by Proposition~\inclusionZEFP{} $(D,\psi)$
would have the ZEFP property.

Support for a positive answer to Question~\AllZEFP{} is provided by
St\o{}rmer's  result~\cite{\StormerZZShiftCstar} that if $D$ is any unital
C$^*$--algebra and $\psi$ is any state on $D$ (with faithful GNS
representation), then letting $(A,\phi)$ be the free product of infinitely many
copies of $(D,\psi)$ indexed
by a set $I$, letting $\sigma$ be a permutation of $I$ without cycles and
letting $\sigma_*$ be the corresponding free permutation of $A$, the
CNT--entropy $h_\phi(\sigma_*)$ of $\sigma_*$ with respect to the free product
state $\phi$ is zero.

\proclaim{Question \AnyFinite}\rm
Given a reduced free product of C$^*$--algebras
$(A,\phi)=(A_1,\phi_1)*(A_2,\phi_2)$, with $\dim(A_1)\ge2$ and
$\dim(A_2)\ge3$ and where $\phi_1$ and $\phi_2$ faithful states, is there an
automorphism $\alpha\in\Aut(A)$, such that $0<ht(\alpha)<\infty\,$?
\endproclaim

It may be especially interesting to restrict the above question to the case
when the states $\phi_1$ and $\phi_2$ are traces.
A first example to consider might be
$(A,\tau)=(C^*_r(\Ints_2),\tau_{\Ints_2})*(C(X),\tau_X)$, where $X$ is the
compact Hausdorff space obtained as the product of infinitely many two--element
spaces and where $\tau_X$ is the state given by the product of uniform
measures.
Now take $\alpha\in\Aut(A)$ to be $\alpha=\id_{C^*_r(\Ints_2)}*\beta$ where
$\beta$ is the Bernoulli shift.
Then $ht(\alpha)\ge ht(\beta)=\log2$.
Is $ht(\alpha)$ finite?

Note, however, that it is easy to find a reduced {\it amalgamated} free product
of C$^*$--algebras $(A,\phi)=(A_1,\phi_1)*(A_2,\phi_2)$ with $A$ non--nuclear
and $\alpha\in\Aut(A)$ with $0<ht(\alpha)<\infty$.
Indeed consider abelian C$^*$--algebras $A_i=C(\Tcirc)\otimes C(X)$, for some
compact Hausdorff space $X$;
let $B=1\otimes C(X)\subseteq A_i$ and let $\phi_i:A_i\to B$ be the slice map
obtained from Haar measure on $\Tcirc$;
then $A=C^*_r(F_2)\otimes C(X)$.
Let $\alpha=\id_{C^*(F_2)}\otimes\beta\in\Aut(A)$, where $\beta$ is an
automorphism on $C(X)$ having strictly positive and finite topological entropy.
By properties of the Brown--Voiculescu topological
entropy~\cite{\BrownZZTopEntropy}, we have
$ht(\beta)\le ht(\alpha)\le ht(\id_{C^*_r(F_2)})+ht(\beta)=ht(\beta)$.

\enddocument